\documentclass[11pt,a4paper]{scrartcl} %{elsarticle}
\KOMAoptions{%
draft=false,
abstract=true,
numbers=enddot
}

\usepackage{amsmath}    %Not needed with amsart
\usepackage{amsthm}     %Not needed with amsart
\usepackage{amsfonts,amssymb}
\usepackage{array}

\usepackage[labelsep=period,justification=centerlast]{caption}%[2011/09/30]
\usepackage{graphicx}

\usepackage[dvipdfmx]{hyperref}
\hypersetup{%
colorlinks=true, linkcolor=blue, citecolor=blue, filecolor=blue, urlcolor=blue,
pdfauthor={A. Kurnosenko}, pdftitle={Lima\c{c}on-like spiral}}

\addtokomafont{section}{\Large}
\addtokomafont{subsection}{\large\itshape}
\DeclareMathOperator{\const}{const}

\newcommand{\Skip}[1]{}

\begin{document}

\newcommand{\kap}{\kappa}
\newsavebox{\TmpBox}
\newcommand{\tmp}{}
\newlength{\tmplength}                    %
%        Labelled/unlabelled equations
%
\newcommand{\Equa}[2]{\begin{equation}#2\label{#1}\end{equation}}
%        NonLabelled equation
\newcommand{\equa}[1]{\[ #1 \]}

\newcommand{\ceq}{\mathrel{\raisebox{0.5pt}{:}\mspace{-3mu}=}}

%        Math. functions
%
\newcommand{\abs}[1]{\left\lvert#1\right\rvert}
\renewcommand{\Re}[1]{\mathop{\rm Re}\nolimits\,(#1)}
\renewcommand{\Im}[1]{\mathop{\rm Im}\nolimits\,(#1)}
\newcommand{\sgn}{\mathop{\rm sgn}\nolimits}
\newcommand{\Deg}[1]{{\ifmmode{#1}^\circ\else${#1}^\circ$\fi}}
\newcommand{\e}{{\mathrm{e}}}
\newcommand{\Exp}[1]{\e^{#1}}
\newcommand{\nvec}[1]{\mathbf{n}(#1)}
%
%       sqrt(-1):
\newcommand{\iu}{{\mathrm{i}}}
\newcommand{\Arc}[1]{\displaystyle{\buildrel\,\,\frown\over{#1}}}
\newcommand{\Brace}[1]{\left\{#1\right\}}
\newcommand{\br}[1]{\left[#1\right]}
\newcommand{\ieq}{\,{=}\,}
\newcommand{\ineq}{\,{\not=}\,}
\newcommand{\ilt}{\,{<}\,}
\newcommand{\igt}{\,{>}\,}
\newcommand{\In}{\,{\in}\,}
\renewcommand{\le}{\leqslant}
\newcommand{\ile}{\,{\le}\,}
\renewcommand{\ge}{\geqslant}
\newcommand{\ige}{\,{\ge}\,}

\newcommand{\Eqref}[1]{\buildrel\eqref{#1}\over=}
%
%            Others
\newcommand{\So}{\quad\Longrightarrow\quad}
\newcommand{\Vect}[1]{\overrightarrow{#1}}
\newcommand{\HM}{\hphantom{{-}}}

%%  Cumulative angles:
%
%\newcommand{\acum}{\widetilde\alpha}
%\newcommand{\bcum}{\widetilde\beta}

%%        Circles [ K(...) ]
%
%\newcommand{\Kls}[1]{{\ifmmode{\cal K}^{\star}_{#1}\else${\cal K}^{\star}_{#1}$\fi}}
\newcommand{\Kl}[1]{{\ifmmode{\cal K}_{#1}\else${\cal K}_{#1}$\fi}}
%\newcommand{\Kr}[1]{K(#1)}

%
%  Refs to figs.
%
\def\FigDir{}
\def\figurename{Fig.}
\newcommand{\Figref}[1]{\ref{F#1}}
\newcommand{\Reffig}[2][]{fig.\;\Figref{#2}\textit{#1}} 
\newcommand{\RefFig}[2][]{Fig.\;\Figref{#2}\textit{#1}} 
\newcommand{\Infigw}[2]{%  Width Id
\includegraphics[width=#1]{\FigDir#2.eps}}

\newcommand{\Infig}[3]{\Infigw{#1}{#2}\caption{{\small #3}}\label{F#2}}

\newcommand{\Bfig}[3]{%  {width}  {id}  (caption)
\parbox[b]{#1}{\Infig{#1}{#2}{#3}}%\HideDisplacementBoxes
}

% from bilens.tex
\newcommand{\Pfig}[3]{%  {width}  {id}  (caption)
%\captionstyle{centerlast}
\centering\Infig{#1}{#2}{#3}}

%==================================================================== MyMath.tex, end

\newcommand{\quo}[1]{\glqq#1\grqq}

\author{Alexey Kurnosenko}
\date{}
%\ead{Alexey.K@dxdy.ru}
%\address{Private research (Russian Federation)}

\title{A nice lima\c{c}on-like spiral}

\maketitle{}%\thispagestyle{plain}

\begin{abstract}
A lima\c{c}on-like curve,
%providing ?
allowing 
$2\pi$-transition with monotone curvature between concentric curvature elements,
is presented.
The curve is 4th degree algebraic, 4th degree rational,
and shares other common features with Pascal's lima\c{c}on.
\end{abstract}\medskip

\noindent
We present a family of curves, traced in polar coordinates $(r,\xi)$ as
\equa{%Polar}{%
    r(\xi)=f\left(\mu\cos\xi+\sqrt{2-\cos^2\xi}\right),\qquad \mu\lesseqgtr0.
}
Dimensionless parameter $\mu$ controls the shape of the curve.
Implicit equation looks like
\Equa{Implicit}{%
     \left(x^2+y^2-\mu f x\right)^2=
     f^2(x^2+ {\color{red}\underline{{\color{black}2}}}y^2),
}
and rational parametrization is
\equa{%Rational}{%
%    t=v\sqrt{\frac{e_0+e_2}{e_0+e_1}}{:}\quad
      x(t)=f\dfrac{\left(t^2-1\right)\,\left[(\mu-1)t^2-(\mu+1)\right]}%
      {t^4+1},\quad
      y(t)=f\sqrt2\, \dfrac{t\left[(\mu-1)t^2-(\mu+1)\right]}{t^4+1}.
}    

Underlined factor $2$ in Eq.\,\eqref{Implicit} distinguishes this curve 
from Pascal's lima\c{c}on \cite{Shikin}, and provides a nice pro\-perty:
{\em extremal circles of curvature are concentric,
and the curve performs $G^3$-continuous transition between them.
If they are equally directed $(\abs{\mu}>1)$, the transition is spiral,
i.\,e. curvature varies monotonically}\footnote{%
{\em Spirality} is meant here in the strong sense, accepted in Computer-Aided Design applications,
assuming monotonicity of curvature.
}.
Total turning angle is in this case~$2\pi$.
%Note that in the case of opposite direction of concentric circles
%spiral transition does not exist.

%\parbox[c]{\textwidth}{\Infigw{\textwidth}{Limacon}}

\begin{figure}[h]
\centering%
\Bfig{.96\textwidth}{LimaconK}%
%{Spiral transition and the plot of curvature~$k$ vs arc length~$s$}
{Spiral transition (heavy curve), its circles of curvature (dashed), their midcircle (dotted-dashed).
Plot of curvature~$k$ vs arc length~$s$.}\end{figure}

The transition is shown as curve $\Arc{AB}$ in \RefFig{LimaconK}.
The whole lima\c{c}on includes the second arc, symmetric about the $x$-axis.
Curvature elements at the endpoints $A$ and~$B$ 
are written below as $\Brace{x,y,\tau,k}$, where $\tau$ defines the unit tangent $(\cos\tau,\sin\tau)^T$,
and~$k$ is curvature.
Common center of two circles is denoted as~$x_c$:
\equa{%CrvEl}{%
  \begin{array}{lllll}
    \xi=0:&  
      \left\{\vphantom{\dfrac{\pi}2}   
         f(\mu+1),\right.& 0, &\dfrac{\pi}2\sgn[f(\mu+1)],&
         \dfrac{2\mu}{\abs{f}\,\abs{\mu+1}(\mu+1)}\left.\vphantom{\dfrac{\pi}2}\right\};\\[1ex]
    \xi=\pi:&  
       \left\{\vphantom{\dfrac{\pi}2}
        f(\mu-1),\right.& 0, &\dfrac{\pi}2\sgn[f(\mu-1)],&
        \dfrac{2\mu}{\abs{f}\,\abs{\mu-1}(\mu-1)}\left.\vphantom{\dfrac{\pi}2}\right\};
  \end{array}
  \qquad   x_c=f\frac{\mu^2-1}{2\mu}.
}
Any ratio of curvatures of concentric circles (except $\pm1$) can be reached with proper~$\mu$:
\equa{%
     \frac{k(\pi)}{k(0)}=\pm\kap^2{:}\quad
     \pm\kap^2=\frac{\mu+1}{\mu-1}\cdot\abs{\frac{\mu+1}{\mu-1}}\So
     \mu=\frac{\abs{\kap}\pm1}{\abs{\kap}\mp1}
}
(the limit case $f\to0$, $\mu\to\infty$, $\mu f=2R=\const$, yields the ratio $+1$,
and lima\c{c}on~\eqref{Implicit} degenerates into duplicated circle $(x^2+y^2-2R x)^2=0$).

The lima\c{c}on is the inverse, with respect to the circle $x^2+y^2=f^2$, 
of conic~\eqref{Conic}, which has excentricity~$e$, focal parameter~$p$,
and the focii at the points $(x_f,0)$:
\Equa{Conic}{%
  2y^2+(1-\mu^2)x^2+2\mu f x -f^2=0,\qquad e=\sqrt{\frac{\mu^2+1}{2}},\quad
  p=\frac12f,\quad
  x_f=\frac{f}{2(\mu\pm e)}.
}
The former vertical axis of symmetry of the conic was equally
a trivial midcircle of two extremal circles of curvature. 
After inversion it appears in \RefFig{LimaconK} as the  midcircle
of two extremal (concentric) circles of curvature,
and as the circle of symmetry of the whole lima\c{c}on.

Note that Pascal's lima\c{c}on was obtained by inversion of conic with the center 
of inversion in the focus. In~\eqref{Conic} the focus is on the $x$-axis
at some distance $x_f\ne0$ from the center of inversion.

The polar equation $\rho(\theta)$ of the lima\c{c}on with 
the pole in the common center~$x_c$ is given by
\equa{
   4\mu^2\rho^2(\theta)-
   4\mu f\left(\cos\theta+\mu\sqrt{\mu^2+\sin^2\theta}\right)\rho(\theta)+
   (\mu^2-1)^2 f^2=0.
}
\Skip{%
\equa{
    2\mu\frac{\rho(\theta)}{f}=  \cos\theta+\mu\sqrt{\mu^2+\sin^2\theta}\pm
    \sin\theta\sqrt{\frac%
    {(\mu^2+1)^2\sin^2\theta+8\mu^2(\mu^2-1)}%
    {(\mu^2-1)\sin^2\theta+2\mu\left(\mu-\cos\theta\sqrt{\mu^2+\sin^2\theta}\right)}}.
}}
\begin{figure}[t]
\centering%
\Bfig{1\textwidth}{Limacon2}{Concentric lima\c{c}ons, produced by inversion of: a)~hyperbola; 
b)~parabola; c,d)~ellipse.}
\end{figure}

%\medskip%
\RefFig{Limacon2} shows the variety of shapes of the lima\c{c}on. 
\begin{itemize}
\item
Curves with $\abs{\mu}>1$ inherit two vertices from the original hyperbola,
as well as spirality of the transition between them.
Point $(0,0)$ is self-intersection.
\item
In the parabolic case $\abs{\mu}=1$ the inner concentric circle degenerates
into a cusp at the point $(0,0)$. Two halves of the lima\c{c}on remain spirals,
as two halves of the parabola were.
\item
When $\abs{\mu}<1$, the concentric circles of curvature have opposite orientations.
Such boundary conditions contradict to spirality.
To enable connection, conic~\eqref{Conic} turns into ellipse,
and each branch of the lima\c{c}on makes use of additional vertex.
Point $(0,0)$ is  isolated singularity.
The special case $\mu=0$
% (with common center $x_c$ at infinity)
is a particluar case of elliptic lemniscate~\cite{Shikin}.
\end{itemize}

The curve was initially obtained from a close look at the critical solution of the method
\cite{AKhyperb}, described there as  $\sigma=\pi$.
Together with concentric given data, the solution 
promised to be simple and interesting.
For normalized (in terms of \cite{AKhyperb}) boundary conditions
\equa{
  \Brace{-1,\: 0,\:{-\dfrac{\pi}2},\:\dfrac{\kap^2-1}{2\kap^2}}\quad\text{and}\quad
  \Brace{ 1,\: 0,\:{-\dfrac{\pi}2},\:\dfrac{\kap^2-1}{2}},\quad
  \text{with common center}\;
  \left(\frac{\kap^2+1}{\kap^2-1},\,0\right),
}
the method returns parametrization
\equa{
  x(u)=\frac{\kap+1}{\kap-1}+
       \frac{2\kap}{\kap-1}\cdot
       \frac{(1-2u)\left[u^2\kap-(1-u)^2\right]}{\kap^2 u^4+(1-u)^4},\quad
  y(u)=\frac{-2\sqrt{2\kap^3}}{\kap-1}\cdot
       \frac{u(1-u)(1-2u)}{\kap^2 u^4+(1-u)^4}.
}
The sought for spiral connection corresponds to $0\le u\le1$.

In common geometric terms, the curve could be constructed as follows.
Consider canonical hyperbola
% (i.\,e., centered in the coordinate origin).
$\frac{x^2}{a^2}-\frac{y^2}{b^2}=1$.
Let~$p$, $e$, and  $z(t)=x(t)+i y(t)$ be its focal parameter, excentricity, and parametrization.
Choose the circle of inversion,
% of radius $R=p$,
centered on the $x$-axis at the point $(x_0,0)=(\mu p,0)$.
The inverse curve $\widetilde{z}(t)$ can be obtained, e.\,g., as
$\frac{p^2}{z(t)-x_0}$,
which also includes reflection about the $x$-axis,
and translation, such that the image of the former infinite point
is shifted from~$(x_0,0)$ to the coordinate origin.
The parameter values~$t_{1,2}$, corresponding to vertices of the hyperbola,
remain such for the curve-image.
Calculating curvature elements at  $\widetilde{z}(t_1)$ and $\widetilde{z}(t_2)$,
and equating the centers of curvature, results in condition 
\equa{2e^2=\mu^2+1.}
So, special choice of excentricity solves the problem of concentricity.


\begin{thebibliography}{9}
\bibitem{Shikin}
Shikin~E.\,V. {\em Handbook and Atlas of Curves.} Boca Raton, FL: CRC Press, 1995.
%
\bibitem{AKhyperb}
Kurnosenko~A.\,I.
\href{http://www.sciencedirect.com/science/article/B6TYN-4YK7HWJ-1/2/7c23f5c3d2d660f495adbde6e606df2a}%
{%
{\em Two-point ${G}^2$ Hermite interpolation with spirals by inversion of hyperbola}.
}
{Comp. Aided Geom. Design}, 27(2010), 474--481.
%
\end{thebibliography}
\end{document}